# Regularity Properties of Constrained Set-Valued Mappings

Marco Papi[*,1,#]   &   Simone Sbaraglia[◇,2,#]

[*]*Istituto per le Applicazioni del Calcolo "M.Picone", V.le del Policlinico 137, I-00161 and Dip. di Matematica Un. "Tor Vergata", Via della Ricerca Scientifica, 00133, Roma, Italy*

[◇]*Istituto per le Applicazioni del Calcolo "M.Picone", V.le del Policlinico 137, I-00161 and Dip. di Metodi e Modelli Matematici Un. "La Sapienza", Via Scarpa 16, 00161, Roma Italy*

**Abstract**

In these notes, we present a general result concerning the Lipschitz regularity of a certain type of set-valued maps often found in constrained optimization and control problems. The class of multifunctions examined in this paper is characterized by means of a set of Lipschitz continuous constraint functions defined on some Lipschitz manifold. The proof of the regularity result for this class of multifunctions is based on a quantitative version of the Implicit Function Theorem for Lipschitzian maps which provides estimates for the neighborhoods where the implicit map can be defined.

## 1   Introduction

This paper deals with the Lipschitz regularity of set-valued maps of the following type:

$$U(x) := \{u \in \mathcal{M} \ : \ c(u,x) \in C\} \tag{1.1}$$

where $x$ lies in a given subset $X$ of $\mathbb{R}^m$, $\mathcal{M}$ is a Lipschitz manifold in $\mathbb{R}^n$ and the function $c = (c_1, \ldots, c_j)$, with the set $C \subset \mathbb{R}^j$, describes some constraint on the variable $u$. In particular (1.1) may include both inequality and equality (or active) constraints.

More precisely, we prove that under compactness assumptions on $\mathcal{M}$ and under independence assumptions on $c$, the map (1.1) is Lipschitz continuous w.r.t. the Hausdorff metric ($d_H$) [4] assigned over the space of all compact subsets of $\mathbb{R}^n$.

[1]e-mail: papi.iac.rm.cnr.it - ivnpp.tiscalinet.it.

[2]e-mail: sbaragli.iac.rm.cnr.it.

[#] This research was supported by the Italian National Research Council (CNR).

The relevance of this result is related to optimization problems, nonlinear analysis, mathematical economics, etc. Lipschitz type properties are important to study mappings associated with the solutions of optimization problems, including variational inequalities, and mathematical programs [13]. The regularity of the value function that arises in stochastic optimization problems with constraints on the admissible controls, depends on the regularity of the *marginal function* [1], [32] and, as is proved in [1], the regularity of the marginal function is connected to the regularity of the multifunction which defines the set of admissible controls, see also [5], [6] and [7].

Our regularity result is also useful to study the Lipschitzian stability of solution maps depending on parameter perturbations. This kind of stability is important not only for a better understanding of the solution behavior w.r.t. perturbations, but also to construct numerical algorithms, see [14], [20], [21], [24].

The Lipschitz continuity is an important analytical property which allows to deal with a variety of problems which arise in set-valued analysis, such as the construction of regular selection maps, for which there is extensive literature elaborating sufficient conditions for the existence of Lipschitz-continuous selections [2], [3], [11]. At the same time, Lipschitz regularity is often too strong and is replaced with weaker conditions such as the upper-Lipschitz continuity, which describes a kind of unilateral Lipschitzian behavior of the multifunctions which does not reduce to the classical local Lipschitz continuity in the case of single-valued maps but turns out to be quite natural in nonlinear optimization; Rockafellar [30] gave a characterization of the upper Lipschitz continuity by employing the so-called proto-derivatives, and [17], [18] contain conditions ensuring the upper Lipschitzian property for a class of multifunctions. Another generalization of the classical Lipschitz-continuity concept, at which many authors paid much attention, is the so-called pseudo-Lipschitzian continuity. This property of multifunctions was introduced by Aubin [1] and used to derive regularity results concerning with solution maps. Various criteria for pseudo-Lipschitz continuity and other related properties of general multifunctions can be found in Mordukhovich [24], [21], [22], [23] and are employed to obtain Lipschitzian stability results in parametric equations. Our approach can be successfully adapted to deal with these more recently introduced regularity concepts.

The typical multifunctions which arise in many practical problems fit into the general format (1.1). When a set is described by equality and inequality constraints with smooth functions, under the well known Mangasarian-Fromovitz constraint qualification [19], [20] ,[28] Lipschitz-type properties can be proved. For instance in [24] a generalized version of that condition is used to derive the pseudo-Lipschitzian behavior of a multifunction such as (1.1), when $\mathcal{M}$ is the whole space $\mathbb{R}^n$ and under some regularity assumptions on the set $C$.

The novelties introduced in our paper lie in the generality of the Manifold $\mathcal{M}$ (that is only Lipschitz continuous), the non smoothness of the functions $c$ and the generality of the set $C$.

Several applications present multifunctions defined by means of inequality constraints for a set of functions defined over a manifold that represents the control space. Unfortunately, in many cases this manifold fails to be differentiable. Actually, the setting which covers a large class of practical situations is the Lipschitzian framework.

The result in this paper extends some previous results proved by the authors in [26], where it is proved the regularity of the optimal value function of a general stochastic control problem in discrete time. In that work, the case of constraints which are $C^1$ on some regions of the Lipschitz manifold $\mathcal{M}$ is treated. Moreover, in [26], is proposed an application of the regularity result to an optimal portfolio allocation problem.

To prove the main result in this paper, we use an extension to manifolds of the generalized notion of Jacobian given by F.H. Clarke [8], [9].

A crucial assumption concerns with the surjectivity of the generalized $u-$Jacobian of the vector valued map $\mathbf{c} = (c_1, \ldots, c_j)$ over $\mathcal{M}$. This condition is quite natural and allows us to use a quantitative formulation of the Implicit Function Theorem that provides an estimate of the neighborhoods where the implicit function is defined.

A version of this Theorem, was stated by Clarke in [8] and [9]. Another, weaker, version of the Theorem was proved in [16] for functions which admit a possible unbounded generalized Jacobian. In both cases the implicit function, which explicits a variable with respect to the other, is locally Lipschitz continuous. Nevertheless these generalizations of the Implicit Function Theorem do not provide a lower estimate of the radius of the neighborhood where the implicit function can be defined. In Theorem 3.4, we prove these kind of estimates (see (3.21)) and as can be seen in the first part of the proof of Theorem 3.2, this bound is essential to prove the Lipschitz continuity of (1.1).

The organization of the paper is as follows. Section 2 contains some backround material about the generalized differentiation of single-valued functions. We recall the main calculus rules which are used in what follows and we introduce a particular function (see (2.8)), which is used in the proof of the Implicit Function Theorem.
In Section 3 we provide the main result concerning the Lipschitz-continuity of (1.1), while Section 4 is devoted to the proofs of the Implicit Function Theorem and of an abstract Lemma.

## 2 Generalized Jacobians

Let us recall the notion of generalized Jacobians for Lipschitz continuous mappings and some fundamental results we need. Here we also generalize the notion of generalized Jacobian to functions whose domain of definition is a Lipschitz manifold.
The ball in $\mathbb{R}^p$ with center $x$ and radius $r$ is denoted by $B_p(x,r)$, in the space $\mathbb{R}^{p \times q}$ of $p \times q$ matrices we use the notation $B_{p \times q}$ for the unit ball with center 0 obtained using the following norm:

$$\|O\| := \max_{\substack{x \in \mathbb{R}^q \\ |x|=1}} |Ox|. \tag{2.2}$$

for every $O \in \mathbb{R}^{p \times q}$.

When a vector-valued map $f : \Omega \to \mathbb{R}^p$, with $\Omega \subset \mathbb{R}^q$ open, is Lipschitz continuous near some point $x \in \Omega$, by Rademacher's Theorem it is differentiable almost everywhere (for the Lebesgue measure) on any neighborhood of $x$ on which $f$ is Lipschitz continuous. We can thus give the following Definition, see [9].

**Definition 2.1** *Let $f : \Omega \to \mathbb{R}^p$, $\Omega \subset \mathbb{R}^q$ open, be Lipschitz continuous near $x \in \Omega$. Then the generalized Jacobian $\partial f(x)$ of $f$ at the point $x \in \Omega$ is defined by*

$$\partial f(x) := co\{P \in \mathbb{R}^{p \times q} : \ P = \lim_{x_h \to x} \frac{\partial f}{\partial x}(x_h), \ f \text{ is differentiable at } x_h \text{ for all } h\},$$

*where the notation "co" indicates the convex hull of a set.*

Some fundamental properties of the generalized Jacobians are summarized below.

**Definition 2.2** *Let $S$ be a set-valued map from $\Omega \subset \mathbb{R}^q$ to $\mathbb{R}^p$. It is said to be upper semicontinuous (u.s.c.) at $x \in \Omega$ if, for every neighborhood $I$ of $S(x) \subset \mathbb{R}^p$, there exists a neighborhood $J$ of $x$ in $\mathbb{R}^q$ such that $S(J \cap \Omega) \subset I$.*

**Proposition 2.3** *([9] Prop. 2.6.2) Let $f : \Omega \to \mathbb{R}^p$, with $\Omega \subset \mathbb{R}^q$ open, be a mapping which is Lipschitz continuous near some $x \in \Omega$. Then the following statements hold:*

- *$\partial f(x)$ is a non-empty convex compact subset of $\mathbb{R}^{p \times q}$.*

- $\partial f(x)$ is closed at $x$.

- $\partial f(x)$ is u.s.c. at $x$. In particular, for any $\varepsilon > 0$, there exists $\delta > 0$ such that for all $y \in B_q(x, \delta) \cap \Omega$
$$\partial f(y) \subset \partial f(x) + \varepsilon B_{p \times q}.$$

We extend Definition 2.1 to maps defined over a Lipschitz manifold $\mathcal{M}$.

**Definition 2.4** *A set $\mathcal{M} \subset \mathbb{R}^n$ is called a d-dimensional **Lipschitz manifold** if for every $u \in \mathcal{M}$ there exists an open neighborhood $W_u$ of $u$ in $\mathcal{M}$ and a homeomorphism $\psi_u : W_u \to V_u$, where $V_u$ is an open subset of $\mathbb{R}^d$ such that $\psi_u$ and $\psi_u^{-1}$ are Lipschitz continuous. The couple $(W_u, \psi_u)$ will be called a chart.*

**Definition 2.5** *An **atlas** on a Lipschitz manifold $\mathcal{M} \subset \mathbb{R}^n$ is a collection of charts $\{(W_\alpha, \psi_\alpha)\}_{\alpha \in A}$ such that $\{W_\alpha\}_{\alpha \in A}$ is a cover of $\mathcal{M}$.*

Hereafter we suppose that every Lipschitz manifold $\mathcal{M}$, is endowed with an atlas, moreover if $\mathcal{M}$ is compact we shall assume, without loss of generality, that the atlas is finite.

**Definition 2.6** *Let $\mathcal{M} \subset \mathbb{R}^n$ be a d-dimensional Lipschitz manifold, and let $f : \mathcal{M} \to \mathbb{R}^p$ be a locally Lipschitz continuous function and let $u \in \mathcal{M}$. For every chart $(W, \psi)$ in the atlas of $\mathcal{M}$, with $u \in W$, let $v_\psi = \psi(u)$ then we define the generalized Jacobian of $f$ at $u$ as*

$$\partial^{\mathcal{M}} f(u) := co \left( \bigcup_{(W,\psi)} \partial(f \circ \psi^{-1})(v_\psi) \right). \tag{2.3}$$

**Remark 2.7** *We observe that Definition 2.6 generalizes Definition 2.1, since $\mathbb{R}^d$ (or every open subset of it) can be viewed as a d-dimensional manifold with a differentiable structure given by the identity map.*
*We omit the dependence on $\mathcal{M}$ of $\partial^{\mathcal{M}} f(u)$, when it is clear from the context the domain of definition of the function $f$.*

Using Definition 2.6, we can formalize a concept of generalized partial derivative for a locally Lipschitz function $f : \mathcal{M}_1 \times \mathcal{M}_2 \to \mathbb{R}^p$, where $M_i \subset \mathbb{R}^{n_i}$ is a Lipschitz manifold whose dimension is $d_i$, $i = 1, 2$. In fact we observe that $\mathcal{M}_1 \times \mathcal{M}_2$ is a Lipschitz manifold whose dimension is $d_1 + d_2$.

**Definition 2.8** *The generalized partial derivative of $f$ at $(u_1, u_2) \in \mathcal{M}_1 \times \mathcal{M}_2$ w.r.t. $u_1$ is defined by the following set*

$$\partial_{u_1} f(u_1, u_2) := \left\{ P \in \mathbb{R}^{p \times n_1} \ : \ \exists \, Q \in \mathbb{R}^{p \times n_2}, \ s.t. \ [P, Q] \in \partial f(u_1, u_2) \right\}. \tag{2.4}$$

We now recall some key results of the theory of Lipschitz continuous mappings.

**Theorem 2.9** *(Mean Value Theorem, [9], Prop. 2.6.5) Let $f : \Omega \to \mathbb{R}^q$ be Lipschitz continuous on an open convex set $\Omega \subset \mathbb{R}^p$, and let $x, y \in \Omega$. Then, there exists a matrix $P \in co\partial f([x, y])$ (where $[x, y]$ stands for the straightline segment connecting $x$ and $y$) such that*

$$f(x) - f(y) = P(x - y). \tag{2.5}$$

The following result about chain rules for generalized Jacobians appears in [9].

**Theorem 2.10** *(Chain rule formula, [9], Cor. 2.6.6) Let $f : \mathbb{R}^q \to \mathbb{R}^p$ be Lipschitz near $x$ and let $g : \mathbb{R}^p \to \mathbb{R}^k$ be Lipschitz continuous near $f(x)$. Then,*

$$\partial(g \circ f)(x) \subset co\{\partial g(f(x)) \circ \partial f(x)\}. \tag{2.6}$$

*If $g$ is differentiable near $f(x)$, then the previous inequality holds with an equality, and the convex hull is superfluous.*

In conclusion of this section, we introduce other notations and definitions which are used in the paper.

- Let $V$ be a set of real square matrices of the same order. We say that $V$ is invertible if every element in $V$ is invertible in the usual sense and we denote by $V^{-1}$ the "inverse" of $V$, that is the set of matrices which are the inverses of elements of $V$.

- Let $V$ be a bounded set of real matrices of the same dimension, then we use the following notation:

$$\|V\|_\infty := \sup_{v \in V} \|v\|, \tag{2.7}$$

where $\|\cdot\|$ is the norm (2.2).

- Let $S$ be a set-valued map defined in an open neighborhood of radius $r$ of $x \in \mathbb{R}^q$, with $S(x) \subset \mathbb{R}^{p \times p}$ and u.s.c. at $x$. Assume that $S(x)$ is compact, convex and invertible with a bounded inverse. Consider the map $\phi(\cdot; S, x) : \delta \mapsto co\{S(B_q(x, \delta))\}$, defined for $\delta \in [0, r)$. This set-valued map is u.s.c. at 0 in the sense of Definition 2.2: in fact, let $I \subset \mathbb{R}^{p \times p}$ be a neighborhood of $\phi(0; S, x)$. Being $S(x)$ convex we have $\phi(0; S, x) = S(x)$ and since $S(x)$ is compact, we find $\varepsilon > 0$ such that $S(x) + \varepsilon B_{p \times p} \subset I$. Since $S$ is u.s.c. at $x$, $S$ is also u.s.c. at $x$ in the $\varepsilon$ sense (see [1], page 45), therefore we find $\eta > 0$ s.t.

$$S(B_q(x, \eta)) \subset S(x) + \varepsilon B_{p \times p}.$$

Hence for every $\delta \in [0, \eta)$, we have

$$\phi(\delta; S, x) = co\{S(B_q(x, \delta))\} \subset \phi(0; S, x) + \varepsilon B_{p \times p} \subset I,$$

where the first inclusion holds because $S(x)$ is supposed to be convex. Let $\lambda > \|S(x)^{-1}\|_\infty$, and

$$\delta(\lambda; S, x) := \sup\left\{\delta > 0 : \|\phi(t; S, \Omega, x)^{-1}\|_\infty \leq \lambda, \ \forall\, t \in [0, \delta)\right\} \tag{2.8}$$

This definition is well-posed, since the set

$$I = I(p, \lambda) := \left\{Q \in \mathbb{R}^{p \times p} : Q^{-1} \text{ exists and } \|Q^{-1}\| < \lambda\right\} \tag{2.9}$$

is open in $\mathbb{R}^{p \times p}$, and $\phi(0; S, x) \subset I$, therefore we can find $\delta > 0$, s.t.

$$\phi(t; S, x) \subset I, \quad \forall\, t \in [0, \delta).$$

This implies that every $Q \in \phi(t; S, x)$ is invertible and $\|Q^{-1}\| < \lambda$. Taking the supremum over $\phi(t; S, x)$, we get $\|\phi(t; S, x)^{-1}\|_\infty \leq \lambda$, for every $0 \leq t < \delta$. This proves that the set on the right-hand side of (2.8) is non empty, and $\delta(\lambda; S, x)$ is well defined.

**Remark 2.11** *It is clear, by the definition of the function (2.8), that the map $S$ is defined in the open ball of centre $x$ and radius $\delta(\lambda; S, x)$.*

- Let $d \geq j$ be two integers and $\Pi(d,j)$ the set of all multi-indexes $\pi = (i_1, \ldots, i_j)$, with $1 \leq i_1 < i_2 < \ldots < i_j \leq d$. For every $\pi = (i_1, \ldots, i_j) \in \Pi(d,j)$, let denote by $T_\pi^{d,j} \in \mathbb{R}^{d \times j}$ the matrix defined by

$$(T_\pi^{d,j})_{h,l} = \begin{cases} 1 & h = i_l \\ 0 & \text{otherwise} \end{cases} \quad (2.10)$$

for every $1 \leq h \leq d$, and $1 \leq l \leq j$. Then for every $\pi \in \Pi(d,j)$ and $H \in \mathbb{R}^{j \times d}$, let denote by $H_\pi$ the $j \times j$ matrix $H \cdot T_\pi^{d,j}$. This matrix is the submatrix of $H$ whose columns correspond to the indexes of $\pi$.

- Given the pair $j, d$, with $j < d$ and $\pi = (i_1, \ldots, i_j) \in \Pi(d,j)$, let $\{\mu_1 < \ldots < \mu_{d-j}\} = \{1, \ldots, d\} \setminus \{i_1, \ldots, i_j\}$ define the following $d-j \times d$ matrix

$$(D_\pi^{d,j})_{h,l} = \begin{cases} 1 & l = \mu_h \\ 0 & \text{otherwise} \end{cases} \quad (2.11)$$

for every $1 \leq h \leq d-j$, and $1 \leq l \leq d$.

## 3 The Lipschitz Regularity

In this section we formulate and prove the main result of the paper concerning the Lipschitz regularity, w.r.t. the Hausdorff metric, of the map (1.1). Let us recall the definition of Hausdorff metric.

**Definition 3.1** *Let $X$ be a separable metric space, with metric $d$. The Hausdorff metric over the space $Comp(X)$ of all compact subsets of $X$ is*

$$d_H(K_1, K_2) = \inf \left\{ \varepsilon > 0 : \ K_1 \subset K_2^\varepsilon, \text{ and } K_2 \subset K_1^\varepsilon \right\} \quad (3.12)$$

*for every $K_1, K_2 \in Comp(X)$. Where for any set $A \subset X$, the set*

$$A^\varepsilon = \{y : \ d(y, A) < \varepsilon\} \quad (3.13)$$

*is the open ball of radius $\varepsilon$ around $A$.*

We refer the interested reader to [4] for a wide treatment of this topic.

**Theorem 3.2** *Let $U$ be the multifunction defined in (1.1), where $\mathcal{M}$ is a compact Lipschitz manifold with dimension $d \geq j$. Make the following assumptions:*

*i) $X$ is a convex subspace of $\mathbb{R}^m$ and $U(x) \neq \emptyset$, for every $x \in X$.*

*ii) The map $c = (c_1, \ldots, c_j)$ is Lipschitz continuous over $\mathcal{M} \times A$, where $X \subset A$, and $A$ is open in $\mathbb{R}^m$.*

*iii) For every $(u, x)$, with $u \in U(x)$ and $x \in X$, there exists $\pi \in \Pi(d,j)$ such that*

$$E(u, x, \pi) := \{P_\pi : \ P \in \partial_u c(u, x)\} \quad (3.14)$$

*is invertible. Furthermore, let*

$$\sup_{\substack{x \in X \\ u \in U(x)}} \sup_{\pi \in \Pi_{u,x}} \|E(u, x, \pi)^{-1}\|_\infty < \tau < \infty, \quad (3.15)$$

where $\Pi_{u,x} := \{\pi : E(u, x, \pi) \text{ is invertible}\}$. Then $x \in X \mapsto U(x)$ is Lipschitz continuous w.r.t. the Hausdorff metric $d_H$ assigned over $Comp(\mathbb{R}^n)$, and its Lipschitz constant can be estimated by a constant $L = L(\tau, Lip(c), Lip_{\mathcal{M}})$, where

$$Lip_{\mathcal{M}} := \max\{Lip(\psi^{-1}) + Lip(\psi) : \psi \text{ is a chart on } \mathcal{M}\}. \tag{3.16}$$

**Remark 3.3** *We note that the amount defined in (3.16) is finite, since the manifold $\mathcal{M}$ is supposed to be compact and therefore the atlas assigned on $\mathcal{M}$ is finite.*

To prove Theorem 3.2 we use the following result which states a quantitative version of the implicit function Theorem for Lipschitz maps providing an explicit estimate of the radii of the balls where the implicit function can be defined.

**Theorem 3.4** *Let $F$ be $\mathbb{R}^d$ valued and Lipschitz continuous on a neighborhood of $(v_0, y_0) \in \mathbb{R}^d \times \mathbb{R}^m$, with $F(v_0, y_0) = 0$. Assume that $\partial_v F(v_0, y_0)$ is invertible with a bounded inverse. Let $Lip(F)$ be the Lipschitz constant of $F$ over the neighborhood of $(v_0, y_0)$. For every*

$$s > 1 + (Lip(F) + 1)\|\partial_v F(v_0, y_0)^{-1}\|_\infty, \tag{3.17}$$

*there exists a constant $\eta_s > 0$ and a continuous function $g : B_m(y_0, \eta_s/2s) \to B_d(v_0, \eta_s)$, such that*

i)

$$g(y_0) = v_0 \tag{3.18}$$

and

$$F(g(y), y) = 0, \qquad \forall\, y \in B_m(y_0, \eta_s/2s). \tag{3.19}$$

*ii) The function $g$ is Lipschitz continuous at $y_0$, i.e.*

$$|g(y) - g(y_0)| \le s|y - y_0|, \quad \forall\, y \in B_m(y_0, \eta_s/2s). \tag{3.20}$$

*iii) The following estimate holds:*

$$\eta_s \ge \frac{1}{2}\delta((s-1)/(Lip(F)+1); \partial_v F, v_0, y_0) \tag{3.21}$$

*where the map $\delta(\cdot\,;\cdot,\cdot)$ is defined in (2.8).*

The proof of Theorem 3.2 relies also on a fundamental abstract Lemma. The proof of Theorem 3.2 is divided into two parts. In the first one, we apply the implicit function Theorem 3.4 to a special function whose structure depends on the relation between the number of contraints ($j$) and the dimension of the manifold ($d$). For every control $u \in U(x)$, we construct a selection map which depends continuously on the state $y$ and we prove that $u$ can be approximated by a control of $U(y)$, for $y$ in a ball centered at $x$. The width of this neighborhood in general may depend on $u$ and $x$. Therefore in the second part of the proof, we use next Lemma to show that if the manifold is compact, the radius of this ball is in fact uniformly bounded in terms of $(u, x)$.

**Lemma 3.5** *Let $\Omega$ be a subset of $\mathbb{R}^q$, and let be given a collection of continuous maps*

$$f_\theta : \Omega_\theta \to Y_\theta \subset \mathbb{R}^k, \quad \theta \in \Theta, \qquad (3.22)$$

*where $\{\Omega_\theta\}_{\theta \in \Theta}$ is a cover of $\Omega$, made of open sets in $\Omega$. For $\theta \in \Theta$ and $\beta \in B$, let $Q_{\theta,\beta}$ be an $\mathbb{R}^\nu$ set-valued map and $S_{\theta,\beta}$ be an $\mathbb{R}^{p \times p}$ set-valued map both defined on $Y_\theta$. Moreover let $J \subset \mathbb{R}^\nu$ be an open set and for $\theta \in \Theta$ and $y \in Y_\theta$ let*

$$V_\theta(y) := \{\beta \in B : \; Q_{\theta,\beta}(y) \subset J\} \qquad (3.23)$$

*be non empty. Make the following assumptions:*

i) *For any $\theta \in \Theta$, $\beta \in B$ and for every $y \in Y_\theta$ such that $\beta \in V_\theta(y)$, $Q_{\theta,\beta}$ is u.s.c. at $y$ and $S_{\theta,\beta}$ is defined in a neighborhood of $y$, in $\mathbb{R}^k$, and is u.s.c. at $y$.*

ii) *For any $\theta \in \Theta$, $\beta \in B$, and $y \in Y_\theta$ with $\beta \in V_\theta(y)$, $S_{\theta,\beta}(y)$ is compact, convex and invertible with a bounded inverse, namely, there exists a constant $\mu > 0$ such that for every $\theta$, $y$ and $\beta \in V_\theta(y)$*

$$\|S_{\theta,\beta}(y)^{-1}\|_\infty < \mu, \qquad (3.24)$$

*Then the function $\Delta : \Omega \to (0, \infty]$, defined by*

$$\Delta(\omega) := \sup_{\theta:\; \omega \in \Omega_\theta} \Delta_\theta(f_\theta(\omega)) \quad \forall\, \omega \in \Omega, \qquad (3.25)$$

$$\qquad (3.26)$$

*where*

$$\Delta_\theta(y) := \sup_{\beta \in V_\theta(y)} \delta(\mu; S_{\theta,\beta}, y) \quad \forall\, y \in Y_\theta,\; \theta \in \Theta, \qquad (3.27)$$

*is lower semicontinuous (l.s.c.).*

We give the proof of Theorem 3.4 and Lemma 3.5 in section 4, while for the remaining part of this section we concentrate on the proof of Theorem 3.2.

**Proof of Theorem 3.2.** Consider $x_1$, $x_2 \in X$. We need to prove that there exists a constant $L > 0$ such that

$$d_H(U(x_1), U(x_2)) \leq L|x_1 - x_2|. \qquad (3.28)$$

Consider the compact straightline segment connecting $x_1$ and $x_2$ and denote it by $K$. The assumption $i)$ says that $K$ is contained in $X$. Now we divide the discussion into two parts.

**First part (construction of a selection map).** Let $u \in U(x)$, with $x \in K$. For every $\pi \in \Pi_{u,x}$ and for every chart $(W, \psi)$ in the atlas assigned on $\mathcal{M}$, such that $u \in W$, we can consider the following Lipschitzian map:

$$F(v, y) := \begin{cases} \begin{bmatrix} c(\psi^{-1}(v), y) - c(u, x) \\ D_\pi^{d,j}(v - \psi(u)) \end{bmatrix} & \text{if } j < d, \\ c(\psi^{-1}(v), y) - c(u, x) & \text{if } j = d, \end{cases} \qquad (3.29)$$

for every $(v, y) \in \psi(W) \times A$. See (2.11) for the definition of the matrix $D_\pi^{d,j}$. We have

$$F(\psi(u), x) = 0. \tag{3.30}$$

Furthermore, since $\mathcal{M} \times A$ can be viewed as a Lipschitz manifold whose charts are given by the "product" of the charts of $\mathcal{M}$ and the identity map over $A$, by definitions (2.3) and (2.4), we have

$$\partial_v F(\psi(u), x) \subset \begin{cases} \begin{bmatrix} \partial_u c(u, x) \\ D_\pi^{d,j} \end{bmatrix} & \text{if } j < d, \\ \partial_u c(u, x) & \text{if } j = d, \end{cases} \tag{3.31}$$

Therefore, if $H \in \partial_v F(\psi(u), x)$, then there exists $P \in \partial_u c(u, x)$ such that

$$H = \begin{cases} \begin{bmatrix} P \\ D_\pi^{d,j} \end{bmatrix} & \text{if } j < d, \\ P & \text{if } j = d, \end{cases} \tag{3.32}$$

By using a development of the determinant of $H$ w.r.t. the set of the first $j$ rows of $H$, by (2.11) it is easy to verify that

$$|\det(H)| = |\det(P_\pi)| > \frac{1}{\tau^j}, \tag{3.33}$$

where the inequality in the right-hand side follows by the assumption $iii)$ and (3.15). So $H$ is invertible. From this observation, we get the invertibility of $\partial_v F(\psi(u), x)$. By (3.33) and since

$$\|\partial_v F(\psi(u), x)\|_\infty \leq Lip(c) \max(Lip_\mathcal{M}, 1) + dj, \tag{3.34}$$

where $Lip_\mathcal{M}$ is the amount defined in (3.16), and using the formula $Cof(H)/det(H)$[3] to determine the inverse of $H$, we obtain that

$$\|\partial_v F(\psi(u), x)^{-1}\|_\infty < \lambda, \tag{3.35}$$

where the constant $\lambda$ depends only on $\tau$, the Lipschitz constant of $c$ and on $Lip_\mathcal{M}$, i.e.

$$\lambda = \lambda(\tau, Lip(c), Lip_\mathcal{M}). \tag{3.36}$$

We can apply Theorem 3.4 to $F$ at $(\psi(u), x)$, choosing the parameter $s$ as

$$s := 1 + \lambda[1 + Lip(c) \max(Lip_\mathcal{M}, 1) + dj], \tag{3.37}$$

In fact by (3.29), it satisfies the inequality (3.17), for every choice of the chart. Therefore there exists a constant $\eta$ and a function $g$ defined on $B_m(x, \eta/2s)$ which takes values in $B_d(\psi(u), \eta)$, such that (3.18) and (3.19) are satisfied. In particular we have

$$c(\psi^{-1}(g(y)), y) = c(u, x) \in C, \quad \forall\, y \in B_m(x, \eta/2s). \tag{3.38}$$

In fact, the inequality follows by $u \in U(x)$. By the inequality (3.38), we get $\psi^{-1}(g(y)) \in U(y)$, for every $y \in X \cap B_m(x, \eta/2s)$. Moreover, by (3.18) and (3.20) we have

$$\begin{aligned}|\psi^{-1}(g(y)) - u| &\leq Lip(\psi^{-1})|g(y) - \psi(u)| = Lip(\psi^{-1})|g(y) - g(x)| \\ &\leq Lip(\psi^{-1})s|y - x| \leq sLip_\mathcal{M}|y - x|. \end{aligned} \tag{3.39}$$

---
[3]The matrix $Cof(\cdot)$ denotes the cofactor matrix.

This previous relation implies that there exists a constant $L(\tau, Lip(c), Lip_{\mathcal{M}})$ such that

$$u \in U(y) + L|x-y|B_n(0,1) = (U(y))^{L|x-y|}, \quad \forall\, y \in X \cap B_m(x, \eta/2s). \tag{3.40}$$

where in the last equality we used the notation (3.13). The same conclusion also holds true for every $(W, \psi)$ such that $u \in W$ and $\pi$ which ensures the inequality (3.35). If we prove that, for a particular choice of $\pi$ and of the chart $(W, \psi)$,

$$\eta > \eta_0 > 0, \tag{3.41}$$

where $\eta_0$ does not depend on $u \in U(x)$ and $x \in K$, we can conclude by (3.40) and Definition 3.1 that $U(\cdot)$ is $d_H$-Lipschitz continuous. The second part of the proof is devoted to proving (3.41).

**Second part (A lower estimate for $\eta$).** We recall that by (3.21) and the particular choice of $s$, it follows

$$\eta \geq \frac{1}{2}\delta(\lambda; \partial_v F, \psi(u), x). \tag{3.42}$$

We shall prove that the number on the righ hand side in (3.42) is bounded from below, by a quantity independent of $(u, x)$. Let

$$\Omega := \{(u, x) : u \in U(x) \text{ and } x \in K\}, \tag{3.43}$$

which is contained in $\mathbb{R}^{n+m}$. Suppose that $\{(W_\theta, \psi_\theta)\}_{\theta \in \Theta}$ is the atlas assigned on $\mathcal{M}$. Set $\Omega_\theta := \Omega \cap (W_\theta \times A)$ and $Y_\theta := \{(\psi_\theta(u), x) : (u, x) \in \Omega_\theta\} \subset \mathbb{R}^{d+m}$, for $\theta \in \Theta$. Since $W_\theta$ is open in $M$ and $\psi_\theta$ is an homeomorphism, $\Omega_\theta$ is open in $\Omega$ and $\{\Omega_\theta\}_{\theta \in \Theta}$ is a cover of $\Omega$. Moreover define $f_\theta : \Omega_\theta \to Y_\theta$ by $f_\theta(u, x) = (\psi_\theta(u), x)$, for every $(u, x) \in \Omega_\theta$, $\theta \in \Theta$. Obviously such maps are continuous as required in Lemma 3.5. Now we proceed with the construction of the collection of multifunctions $S_{\cdot,\cdot}$ and $Q_{\cdot,\cdot}$. Consider $B := \Pi(d, j)$ and let $\theta \in \Theta$, then set

$$\varphi_\theta(v, y) := c(\psi_\theta^{-1}(v), y), \quad \forall\, (v, y) \in \psi(W_\theta) \times A \tag{3.44}$$

which is Lipschitz continuous on its domain of definition, and for $\beta \in B$ let

$$Q_{\theta,\beta}(v, y) := \partial_v \varphi_\theta(v, y) \cdot T_\beta^{d,j} \tag{3.45}$$

for every $(v, y) \in Y_\theta$. Let $J := I(j, \tau)$, where $I(\cdot, \cdot)$ is defined in (2.9). Let $(\theta, \beta) \in \Theta \times B$, then

$$S_{\theta,\beta}(v, y) := \begin{cases} \begin{bmatrix} \partial_v \varphi_\theta(v, y) \\ D_\beta^{d,j} \end{bmatrix} & \text{if } j < d, \\ \partial_v \varphi_\theta(v, y) & \text{if } j = d, \end{cases} \tag{3.46}$$

for every $(v, y) \in \psi_\theta(W_\theta) \times A$. With these definitions the set $V_\theta(v, y)$ is non empty, whenever $(v, y) \in Y_\theta$, $\theta \in \Theta$. In fact, by Definition 2.3 and (3.14), if $(v, y) \in Y_\theta$, it hold $\psi_\theta^{-1}(v) \in U(y)$ and

$$Q_{\theta,\beta'}(v, y) \subset E(\psi_\theta^{-1}(v), y, \beta'), \tag{3.47}$$

which by assumption $iii)$ is invertible for some $\beta' \in \Pi_{\psi_\theta^{-1}(v), y}$. Furthermore by (3.15), the norm of the inverse of every element of $Q_{\theta,\beta'}(v, y)$ is strictly less than $\tau$. So $\beta' \in V_\theta(v, y)$.

Since $\varphi_\theta$ is Lipschitz continuous over $\psi_\theta(W_\theta) \times A$, by Proposition 2.3, $S_{\theta,\beta}$ and $Q_{\theta,\beta}$ are both u.s.c at every point of $Y_\theta$, and $S_{\theta,\beta}$ is defined in an open set which contains $Y_\theta$, for every $\theta$ and $\beta$. This proves that condition $i)$ of Lemma 3.5 is satisfied. To verify the assumption $ii)$ of Lemma 3.5, consider $(v, y) \in Y_\theta$, then the definition of the generalized Jabobian says that $S_{\theta,\beta}(v, y)$ is compact and convex, for every $\beta$.

If $\beta \in V_\theta(v, y)$, every matrix $H \in S_{\theta,\beta}(v, y)$ is written as in (3.32) where $\pi = \beta$ and $P \in \partial_v \varphi_\theta(v, y)$. Since $P_\beta = P \cdot T_\beta^{d,j} \in Q_{\theta,\beta}(v, y)$, proceeding as in the first part of the proof, we infer that $H$ is invertible. Moreover the norm of the inverse of $S_{\theta,\beta}(v, y)$ is strictly less than $\lambda$, where $\lambda$ is defined in (3.36). This yields the inequality (3.24), with $\mu = \lambda$.

Applying Lemma 3.5, the function $\Delta$ defined in (3.25) is l.s.c, therefore it admits a minimum $\Delta_0 > 0$, over $\Omega$: this follows by the compactness of $\Omega$ which can be easily derived by the compactness of the manifold $\mathcal{M}$. For every $x \in K$, $u \in U(x)$, there exists $(W_\theta, \psi_\theta)$, such that $u \in W_\theta$, and some $\beta \in V_\theta(v, x)$, where $v := \psi_\theta(u)$, such that

$$\delta(\lambda; S_{\beta,\pi}, f_\theta(v, x)) > \frac{\Delta_0}{2}. \tag{3.48}$$

Choosing $\psi = \psi_\theta$ and $\pi = \beta$, in the definition (3.29) of $F$, we satisfy the inequality (3.35) and we get

$$\partial_v F(\psi(u), x) = S_{\theta,\beta}(v, x), \tag{3.49}$$

which implies, with (3.48),

$$\delta(\lambda; \partial_v F, \psi(u), x) > \frac{\Delta_0}{2}. \tag{3.50}$$

By (3.42), we have proved the assertion (3.41) with $\eta_0 = \Delta_0/4$.

**Conclusions.** We have shown that

$$\tilde{U}(s) \subset \tilde{U}(t) + L|x_2 - x_1||s - t|, \quad \forall\, s, t \in [0, 1], \text{ s.t. } |x_2 - x_1||t - s| < \frac{\eta_0}{2s}, \tag{3.51}$$

where $\tilde{U}(t) := U(x(t))$, and $x(t) = x_1 + t(x_2 - x_1)$, for $t \in [0, 1]$. Let $h > \frac{2s|x_2-x_1|}{\eta_0}$ be an integer and define $t_i = \frac{i}{h}$, for $i = 0, \ldots, h$. Then (3.51) yields

$$U(x_1) = \tilde{U}(0) \subset \tilde{U}(t_1) + L\frac{|x_2 - x_1|}{h} \subset \cdots \subset \tilde{U}(t_i) + iL\frac{|x_2 - x_1|}{h}$$

$$\subset \tilde{U}(1) + L|x_2 - x_1| = U(x_2) + L|x_2 - x_1|. \tag{3.52}$$

This implies the inequality (3.28) and the proof is complete.

∎

**Remark 3.6** *The relation (3.51) and the chain of inclusions (3.52) imply that the assumption i) of Theorem 3.2, can be weakened by supposing that the space $X$ is connected in $\mathbb{R}^m$ and for every pair of distinct points $x$ and $y$ in $X$, there exists a continuous arc $\gamma : [a, b] \to X$, such that $\gamma(a) = x$, $\gamma(b) = y$ and $\gamma$ is piecewise differentiable with a bounded derivative on the interval $[a, b]$. Furthermore the following holds:*

$$l(\gamma) := \int_a^b |\gamma(t)|dt \leq a(X)|x - y|. \tag{3.53}$$

*Here $a(x) > 0$ is a characteristic constant of the space $X$, which does not depend on the particular arc $\gamma$.*

*The inequality (3.53) is for instance satisfied if $X$ is a compact connected differentiable manifold in $\mathbb{R}^m$.*

**Remark 3.7** *We obesrve that also in presence of a Lipschitz Manifold with a boundary whose dimension d satisfies*

$$d \geq j+1 \tag{3.54}$$

*the same technique used to prove Theorem 3.2 can be successfully applied. Actually, when the point $u \in U(x)$ belongs to the boundary of $\mathcal{M}$, we have a Lipschitz chart on $u$, whose restriction to the boundary is a Lipschitz homeomorphism onto an open set of $\mathbb{R}^{d-1}$ and its inverse is also Lipschitz continuous. Therefore, under the condition (3.54), we can apply exactly the same argument developed in the first part of the proof of Theorem 3.2.*

## 4 The Implicit Function Theorem

In this section we present the proof of the quantitative version of the implicit function Theorem 3.4. It is based on the following Proposition which represents a specialization of Theorem 5.2 in [16]. In that Theorem it is proved that there exist constants $\delta > 0$ and $\varepsilon > 0$ such that the inequality (4.55) holds true when $\frac{1}{\lambda}$ is replaced by $\varepsilon$. In our formulation we prove that by reinforcing the invertibility assumption on the generalized Jacobian, we can take in fact $\frac{1}{\lambda}$ as the parameter $\varepsilon$ and also we show how to choose the parameter $\delta$ as a function of $\lambda$.

**Proposition 4.1** *Let $f$ be a Lipschitz continuous map on a neighborhood of $\xi_0 \in \mathbb{R}^p$, which takes values in $\mathbb{R}^p$. Suppose that $\partial f(\xi_0)$ is invertible with a bounded inverse. For every $\lambda > \|\partial f(\xi_0)^{-1}\|_\infty$, let $\delta_\lambda = \delta(\lambda; \partial f, \xi_0)$ (see (2.8)), then the following holds:*

$$|f(\xi_0 + h) - f(\xi_0)| \geq \frac{1}{\lambda}|h|, \quad \forall h \neq 0, \; |h| < \delta_\lambda \tag{4.55}$$

*and*

$$f(\xi_0) + \frac{\delta}{2\lambda} B_p(0,1) \subset f(\xi_0 + \delta B_p(0,1)), \tag{4.56}$$

*for every $\delta < \delta_\lambda$.*

**Proof of Proposition 4.1.** To prove the first part (4.55), consider $\lambda > \|\partial f(\xi_0)^{-1}\|_\infty$, then consider $h \neq 0$ and $|h| < \delta_\lambda$. By Remark 2.11, $f$ is Lipschitz continuous over $B_p(\xi_0, \delta_\lambda)$, and by the Mean Value Theorem 2.9, we can write

$$f(\xi_0 + h) - f(\xi_0) = Qh, \tag{4.57}$$

where $Q \in co\{\partial f([\xi_0, \xi_0 + h])\}$. By

$$co\{\partial f([\xi_0, \xi_0 + h])\} \subset co\{\partial f(B_p(\xi_0, |h|))\} = \phi(|h|; \partial f, \xi_0) \tag{4.58}$$

and (2.8), we deduce that $Q$ is invertible and $\lambda \geq \|Q^{-1}\|$. Therefore

$$\frac{|Qh|}{|h|} \geq \frac{1}{\|Q^{-1}\|} \geq \frac{1}{\lambda}. \tag{4.59}$$

This proves the inequality (4.55).

The remaining part of the proof follows by applying inequality (4.55), Theorem 2.10 and using exactly the same argument used for proving Theorem 5.2 in [16].

■

We shall denote by $I_m$ the identity matrix over $\mathbb{R}^m$ and by $O_{m,d}$ the $m \times d$ null matrix.

**Proof of Theorem 3.4**. Let $F$ be Lipschitz continuous over $\Omega_0 := B_{d+m}((v_0, y_0), r_0)$, with a Lipschitz constant $Lip(F)$. Consider the function $f(v, y) = (F(v, y), y)$, then

$$\partial f(v, y) = \begin{pmatrix} \partial F(v, y) \\ O_{m,d} \quad I_m \end{pmatrix}$$

for every $(v, y) \in \Omega_0$. If $H \in \partial f(v_0, y_0)$, we can write

$$H = \begin{pmatrix} P & Q \\ O_{m,d} & I_m \end{pmatrix}$$

for some $[P, Q] \in \partial F(v_0, y_0)$. This implies $P \in \partial_v F(v_0, y_0)$ and therefore $H$ is invertible with inverse

$$H^{-1} = \begin{pmatrix} P^{-1} & -P^{-1}Q \\ O_{m,d} & I_m \end{pmatrix}.$$

Since $\|Q\| \leq Lip(F)$, we get

$$\|H^{-1}\| \leq 1 + (1 + Lip(F))\|\partial_v F(v_0, y_0)^{-1}\|_\infty, \tag{4.60}$$

Therefore $\partial f(v_0, y_0)$ is invertible with a bounded inverse and we can apply Proposition 4.1 to $f$ on $\Omega_0$ and at the point $\xi_0 = (v_0, y_0)$. Let $s$ be chosen as in (3.17), therefore by (4.60), the inclusion (4.56) holds with $\delta = \eta_s := \frac{1}{2}\delta_s$. We define the map $g$ in the following way: let $y \in B_m(y_0, \frac{\eta_s}{2s})$, then

$$(0, y) \in f(v_0, y_0) + B_{d+m}((0, 0), \frac{\eta_s}{2s})$$

in fact, $f(v_0, y_0) = (0, y_0)$. By (4.56), we find $(z, w) \in (v_0, y_0) + B_{d+m}((0, 0), \eta_s)$, such that

$$f(z, w) = (0, y). \tag{4.61}$$

We set $g(y) = z$. We have $w = y$, and $|z - v_0| \leq |(z - v_0, w - y_0)| < \eta_s$, therefore $z \in B_d(v_0, \eta_s)$, moreover $z$ is uniquely defined. In fact, if $z, z'$, satisfy $(z, y), (z', y) \in (v_0, y_0) + B_{d+m}((0, 0), \eta_s)$ and both give (4.61), then by the Mean Value Theorem 2.9 applied to $f$, we have

$$0 = f(z', y) - f(z, y) = Q \cdot \begin{pmatrix} z' - z \\ 0 \end{pmatrix} \tag{4.62}$$

where $Q \in co\{\partial f(B_{d+m}(\xi_0, \eta_s))\} = \phi(\eta_s; \partial f, \xi_0)$. Since $\eta_s < \delta(s; \partial f, \xi_0)$, by (2.8) we deduce that $Q$ is invertible. Hence by (4.62) we have $z = z'$. Hence $g$ is well defined over $B_m(y_0, \eta_s/2s)$ and takes values in $B_d(v_0, \eta_s)$. (4.61) and the definition of $f$ yield (3.18) and (3.19), while, by (4.55) we have

$$|y - y_0| = |(0, y) - (0, y_0)| = |f(g(y), y) - f(g(y_0), y_0)| \geq \frac{1}{s}|g(y) - g(y_0)|, \tag{4.63}$$

for every $y \in B_m(y_0, \eta_s/2s)$. We prove now the last assertion (3.21). Let $r_0 > \delta > 0$ such that

$$\|\phi(t; \partial_v F, v_0, y_0)^{-1}\|_\infty \leq \frac{s-1}{1 + Lip(F)}$$

for every $t \in [0, \delta)$. We have

$$co\{\partial f(B_{d+m}(\xi_0, t))\} = \begin{pmatrix} co\{\partial F(B_{d+m}((v_0, y_0), t))\} \\ O_{m,d} \quad I_m \end{pmatrix} \quad (4.64)$$

and observing that

$$co\{\partial_v F(B_{d+m}((v_0, y_0), t))\} = \left\{ P \in \mathbb{R}^{d \times d} : \exists Q \text{ s.t. } [P, Q] \in co\{\partial F(B_{d+m}((v_0, y_0), t))\} \right\},$$

we can repeat the arguments used in the beginning of the proof and by (4.60) and (4.64) deduce that

$$\|\phi(t; \partial f, \xi_0)^{-1}\| \leq s$$

for every $t \in [0, \delta)$. By the definition (2.8), this implies

$$\eta_s = \frac{1}{2}\delta(s; \partial f, \xi_0) \geq \frac{1}{2}\delta((s-1)/(Lip(F)+1); \partial_v F, v_0, y_0).$$

∎

**Proof of Lemma 3.5**. Since the function (3.25) is strictly positive, to prove the assertion it suffices to show that for every $\alpha > 0$, the set $\{\omega \in \Omega : \Delta(\omega) > \alpha\}$ is open in $\Omega$. To this end, we prove this property at first for the collection $\{\{y \in Y_\theta : \Delta_\theta(y) > \alpha\}\}_{\theta \in \Theta}$. Let $\theta \in \Theta$ and $y \in Y_\theta$ be such that $\Delta_\theta(y) > \alpha$. We prove that $y$ admits an open neighborhood in $Y_\theta$ where this inequality holds. By the definition of $\Delta_\theta$, there exists $\beta' \in V_\theta(y)$ such that

$$\delta_\mu := \delta(\mu; S_{\theta,\beta'}, y) > \alpha. \quad (4.65)$$

By the assumptions, $Q_{\theta,\beta'}$ is u.s.c. at $y$ for such $\beta'$, therefore since $J \subset \mathbb{R}^\nu$ is open, by (3.23) we find $r > 0$ such that

$$Q_{\theta,\beta'}(x) \subset J, \quad \forall\, x \in B_k(y, r) \cap Y_\theta, \quad (4.66)$$

Fix $\mu' > 0$ so that $\|S_{\theta,\beta'}(y)^{-1}\|_\infty < \mu' < \mu$. Since $\beta' \in V_\theta(y)$, by the assumption $i)$, we may assume that $S_{\theta,\beta'}$ is defined in an open ball of radius $r' < r$ and centre $y$; moreover this function is u.s.c. at $y$, therefore we can find $\rho \in (0, r')$ such that

$$S_{\theta,\beta'}(x) \subset I(p, \mu'), \quad \forall\, x \in B_k(y, \rho). \quad (4.67)$$

See definition (2.9) for the description of $I(\cdot, \cdot)$. So whenever $x \in B_k(y, \rho) \cap Y_\theta$, (4.66) and (4.67) hold true. In particular $\beta' \in V_\theta(x)$ and so $S_{\theta,\beta'}$ is defined in a neighborhood of $x$, is u.s.c. at $x$, and by $ii)$ $S_{\theta,\beta'}(x)$ is compact, convex and invertible with a bounded inverse and we can evaluate the function (2.8) at $(\mu; S_{\theta,\beta'}, x)$. In fact by (4.67) it follows

$$\|S_{\theta,\beta'}(x)^{-1}\|_\infty \leq \mu' < \mu.$$

Let $\eta > \alpha$ and $\rho > \sigma > 0$ be such that $\delta_\mu > \eta + \sigma$. Then for every $t \in [0, \eta)$ and $x \in B_k(y, \sigma) \cap Y_\theta$ we have $B_k(x, t) \subset B_k(y, \eta + \sigma)$, which implies

$$\phi(t; S_{\theta,\beta'}, x) = co\{S_{\theta,\beta'}(B_k(x, t))\} \subset co\{S_{\theta,\beta'}(B_k(y, \eta + \sigma))\} = $$
$$= \phi(\eta + \sigma; S_{\theta,\beta'}, y),$$

Since $\eta + \sigma < \delta_\mu$, we have

$$\|\phi(t; S_{\theta,\beta'}, x)^{-1}\|_\infty \leq \mu, \quad \forall\, t \in [0, \eta), \tag{4.68}$$

and by this inequality, we obtain

$$\delta(\mu; S_{\theta,\beta'}, x) \geq \eta > \alpha, \quad \forall\, x \in B_k(y, \sigma) \cap Y_\theta. \tag{4.69}$$

By (4.69), whenever $x$ belongs to $B_k(y, \sigma) \cap Y_\theta$,

$$\Delta_\theta(x) = \sup_{\beta \in V_\theta(x)} \delta(\mu; S_{\theta,\beta}, x) \geq \delta(\mu; S_{\theta,\beta'}, x) > \alpha. \tag{4.70}$$

Therefore $y$ is an interior point for the set $\{y \in Y_\theta : \Delta_\theta(y) > \alpha\}$. By the arbitrary choice of $y$, and $\theta$, we have proved the introductory assertion on the collection of sets. Now we proceed to get the same conclusion on $\Delta$.

To simplify the presentation, for every $\omega \in \Omega$, we denote with $\Theta(\omega)$ the set whose elements are the indexes $\theta \in \Theta$ which satisfy $\omega \in \Omega_\theta$.

Let $\omega \in \Omega$, such that $\Delta(\omega) > \alpha$, then we find $\theta' \in \Theta(\omega)$ such that $\Delta_{\theta'}(f_{\theta'}(\omega)) > \alpha$. The previous discussion implies the existence of an open neighborhood $N$ of $f_{\theta'}(\omega)$ contained in $Y_{\theta'}$ where $\Delta_{\theta'} > \alpha$. Consider

$$N' := \{\omega' \in \Omega_{\theta'} : f_{\theta'}(\omega') \in N\},$$

Obviously $\omega \in N'$ and $N'$ is open in $\Omega_{\theta'}$, because of the continuity of $f_{\theta'}$, then it is also open in $\Omega$: in fact $\Omega_{\theta'}$ is supposed to be open in $\Omega$. Then we have

$$\Delta(\omega') = \sup_{\theta \in \Theta(\omega')} \Delta_\theta(f_\theta(\omega')) \geq \Delta_{\theta'}(f_{\theta'}(\omega')) > \alpha, \quad \forall\, \omega' \in N'. \tag{4.71}$$

Actually for every $\omega' \in N'$, $\theta' \in \Theta(\omega')$. The relation (4.71) shows that $\omega$ is an interior point of $\{\omega \in \Omega : \Delta(\omega) > \alpha\}$. Being $\omega$ arbitrary, $\Delta$ is l.s.c. on $\Omega$.

∎

# References

[1] J.P. Aubin, Lipschitz behavior of solutions to convex minimization problems, Math. Oper. Res. 9 (1984) 87-111.

[2] J.P. Aubin, A. Cellina, Differential Inclusions, Set-Valued Maps and Viability Theory, Grundlehren der Mathematischen Wissenschaften, 264, Springer Verlag, Berlin, 1984.

[3] J.P. Aubin, H. Frankowska, Set-Valued Analysis, Birkhäuser, Boston 1990.

[4] M. F. Barnsley, Fractals Everywhere, Second Edition, Academic Press Professional, 1993.

[5] D. Bertsekas, Dynamic Programming, deterministic and stochastic models, Prentice-Hall Inc., 1987.

[6] D. Bertsekas, S.E. Shreve, Stochastic Optimal Control: The Discrete Time Case, Academic Press Inc., 1978.

[7] J.R. Birge, F. Louveaux, Introduction to Stochastic Programming, Springer Series Operations Research, Springer, 1999.

[8] F.H. Clarke, On the inverse function Theorem, Pacific J. Math. 64 (1976) 97-102.

[9] F.H. Clarke, Optimization and non-smooth analysis, Classics in Applied Mathematics, 5, SIAM, Philadelphia, PA, 1990.

[10] G. Debreu, Theory of Value, Wiley, 1959.

[11] G. Dommisch, On the existence of Lipschitz-continuous and differentiable selections for multifunctions, Math. Res. 35 (1987) 60-73.

[12] A.L. Dontchev, W.W. Hager, A new approach to Lipschitz continuity in state constrained optimal control, Systems Control Lett. 35 (1998) 137-143.

[13] A.L. Dontchev, W.W. Hager, Lipschitzian Stability for State Constrained Nonlinear Optimal Control, SIAM J. Control Optim. 36 (1998) 698-718.

[14] P.T. Harker, J.S. Pang, Finite-dimensional variational inequality and nonlinear complementarity problems: a survey of theory, algorithms and applications, Math. Program. 48 (1990) 161-220.

[15] V. Jeyakumar, D.T. Luc, Nonsmooth calculus, minimality and monotonicity of convexificators, J. Optim. Theory Appl. 101 (1999) 599-621.

[16] V. Jeyakumar, D.T. Luc, An Open Mapping Theorem using unbounded Generalized Jacobians, Nonlinear Anal. 50 (2002) 647-663.

[17] A.J. King, R.T. Rockafellar, Sensitivity analysis for nonsmooth generalized equations, Math. Program. 55 (1992) Ser. A 193-212.

[18] A.B. Levy, Implicit multifunction theorems for the sensitivity analysis of variational conditions, Math. Program. 74 (1996) Ser. A 333-350.

[19] O.L. Mangasarian and S. Fromovitz, The Fritz-John necessary optimality conditions in the presence of equality and inequality constraints, J. Math. Anal. Appl. 17 (1967) 37-47.

[20] B. Mordukhovich, Approximation methods in problems of optimization and control, Nauka, Moscow, 1988.

[21] B. Mordukhovich, Sensitivity analysis in nonsmooth optmization, Theoretical Aspects of Industrial Design, SIAM Proceedings in Applied Math. 58 (1992) 32-46.

[22] B. Mordukhovich, Complete characterization of openess, metric regularity, and Lipschitzian properties of multifunctions, Trans. Amer. Math. Soc. 340 (1993) 1-35.

[23] B. Mordukhovich, Lipschitzian stability of constraint systems and generalized equations, Nonlinear Anal. 22 (1994) 173-206.

[24] B. Mordukhovich, Stability analysis for parametric constraint and variational systems by means of set-valued differentiation, Optimization 31 (1994) 13-46.

[25] M. Papi, S. Sbaraglia, Optimal Asset-Liability Management with Constraints: A Dynamic Programming Approach, Quaderno IAC, 2001, submitted to Quantitative Finance.

[26] M. Papi, S. Sbaraglia, Lipschitzian Estimates in Discrete-Time Stochastic Optimal Control, Quaderno IAC, 2002, http://front.math.ucdavis.edu/math.OC/0206037, submitted to SIAM J. Control Optim.

[27] S.M. Robinson, An implicit-function theorem for a class of nonsmooth functions, Math. Oper. Res. 16 (1991) 292-309.

[28] R.T. Rockafellar, The thoery of subgradients and its applications to optimization: convex and nonconvex functions, Heldermann, 1981.

[29] R.T. Rockafellar, Lipschitz properties of multifunctions, Nonlinear Anal. TMA 9 (1985) 867-885.

[30] R.T. Rockafellar, Proto-differentiability of set-valued mappings and its applications in optimization, Ann. Inst. H. Poincoré Anal. Non Linéaire 6 (1989) suppl. 449-482.

[31] R.T. Rockafellar, R. Wets, Variational Analysis, Springer Verlag, Berlin, 1997.

[32] R. Vinter, Optimal Control, Birkhäuser, Boston, 2000.